\documentclass{amsart}
\usepackage{amssymb,graphicx}
\title[Rigidity of bi-Lipschitz equivalence of function-germs]
{Rigidity of bi-Lipschitz equivalence of weighted homogeneous function-germs in the plane}

\author[Alexandre Fernandes]{Alexandre Fernandes}
\address{Departamento de Matem\'atica, Universidade Federal
do Cear\'a, Av. Mister Hull s/n,Campus do PICI, Bloco 914,CEP:
60.455-760 - Fortaleza - CE - Brasil.}
\email{alexandre.fernandes@ufc.br}

\author[Maria Ruas]{Maria Ruas}
\address{Instituto de Ci\^encias Matem\'aticas e Computa\c c\~ao, Av. Trabalhador S\~ao-carlense 400, Centro
Caixa Postal: 668 CEP 13560-970, S\~ao Carlos SP, Brasil}
\email{maasruas@icmc.usp.br}
\date{\today}



\def \N {\mathbb{N}}

\def \C {\mathbb{C}}
\def \Fx {\frac{\partial F}{\partial x}}
\def \Fy {\frac{\partial F}{\partial y}}
\def \Ft {\frac{\partial F}{\partial t}}
\def \Gx {\frac{\partial G}{\partial x}}
\def \Gy {\frac{\partial G}{\partial y}}
\def \Gs {\frac{\partial G}{\partial s}}

\newtheorem{theorem}{Theorem}[section]

\newtheorem{proposition}[theorem]{Proposition}

\theoremstyle{definition}

\newtheorem{example}[theorem]{Example}
\theoremstyle{remark}

\numberwithin{equation}{section}

\begin{document}

\maketitle

\begin{abstract} The main goal of this work is to show that if two weighted homogeneous (but not homogeneous) function-germs $(\C^2,0)\rightarrow(\C,0)$ are bi-Lipschitz equivalent, in the sense that these function-germs can be included in a strongly bi-Lipschitz trivial family of weighted homogeneous function-germs, then they are analytically equivalent.
\end{abstract}

\section{Introduction}

$$f_t(x,y)=xy(x-y)(x-ty) \ ; \ 0<|t|<1$$ defines a family of function-germs $(\C^2,0)\rightarrow(\C,0)$ with isolated singularity. In 1965, H. Whitney justified the rigidity of the analytic classification of function-germs by proving that: given $t\neq s$ there is no $\phi\colon(\C^2,0)\rightarrow(\C^2,0)$ germ of bi-analytic map such that $f_t=f_s\circ\phi$, i.e. $f_t$ is not analytically equivalent to $f_s$. In another way, with respect to the topological point of view, this family is not so interesting, since for any $t\neq s$ there exists a germ of homeomorphism $\phi\colon(\C^2,0)\rightarrow(\C^2,0)$ such that $f_t=f_s\circ\phi$, i.e. $f_t$ is topologically equivalent to $f_s$. In fact, the topological classification of
reduced polynomial function-germs $(\C^2,0)\rightarrow(\C,0)$ is well-understood, as it was shown in \cite{P}.
In the seminal paper \cite{HP}, Henry and Parusinski considered the bi-Lipschitz equivalence, which is between the analytic and the topological equivalence, of function-germs. This paper motivated other papers about the problem of bi-Lipschitz classification of function-germs. For instance,   \cite{PK} and \cite{BFP} showed that, in some sense, for weighted homogeneous real function-germs in two variables the problem of bi-Lipschitz classification is quite close to the problem of analytic classification. The results presented in \cite{HP} point out to a rigidity of the bi-Lipschitz classification of function-germs.  More precisely, they considered the family
$$f_t(x,y)=x^3+y^6-3t^2xy^4 \ ; \ 0<|t|<\frac{1}{2}$$ and proved that: given $t\neq s$, there is no $\phi\colon(\C^2,0)\rightarrow(\C^2,0)$ germ of bi-Lipschitz map such that $f_t=f_s\circ\phi$, i.e. $f_t$ is not bi-Lipschitz equivalent to $f_s$.  The strategy used by them was to introduce a new invariant based on the observation that the bi-Lipschitz homeomorphism does not move much the regions around the relative polar curves.  For a single germ $f$ the invariant  is given in terms of the leading coefficients of the asymptotic expansions of $f$ along the branches of its generic polar curve. In the case that the bi-Lipschitz triviality of a family of function-germs comes by integrating a Lipschitz vector field, here called strong bi-Lipschitz triviality, the calculations are much easier, and very illustrative.
The main goal of this work is to show that if two weighted homogeneous (but not homogeneous) function-germs $(\C^2,0)\rightarrow(\C,0)$ are bi-Lipschitz equivalent, in the  sense that these function-germs can be included in a strongly bi-Lipschitz trivial family of weighted homogeneous function-germs, then they are analytically equivalent.

\section{Preliminaries}

\subsection{Analytic and bi-Lipschitz equivalences}
Two function-germs $f,g\colon (\C^n,0)\rightarrow (\C,0)$ are called \emph{analytically equivalent} if there exists a a germ of bi-analytic map $\phi\colon (\C^n,0)\rightarrow (\C^n,0)$ such that $f=g\circ\phi$.

Two function-germs $f,g\colon (\C^n,0)\rightarrow (\C,0)$ are called \emph{bi-Lipschitz equivalent} if there exists a bi-Lipschitz map-germ $\phi\colon (\C^n,0)\rightarrow (\C^n,0)$ such that $f=g\circ\phi$.

Let $f_t$ ($t\in U$ domain in $\C$) be a family of analytic function-germs. That is, there is a neighborhood
$W$ of $0$ in $\C^n$ and an analytic function $F\colon W\times U\rightarrow\C$ such that $F(0,t)=0$ and $f_t(x)=F(x,t)$ $\forall \ t\in U$, $\forall \ x\in W$. We call $f_t$ \emph{strongly bi-Lipschitz trivial} when there is a continuous family of Lipschitz vector fields $v_t\colon\C^n\rightarrow\C^n$ such that
$$\frac{\partial f_t}{\partial t}(x)=(df_t)(x)(v_t(x))$$
$\forall \ t\in U$, $\forall \ x\in W$.

The following result shows that if we include two function-germs $f$ and $g$ into a strongly bi-Lipschitz trivial family of function-germs, then the two initial function-germs $f$ and $g$ are bi-Lipschitz equivalent.

\begin{theorem}\label{thom-levine} If $f_t$ is strongly bi-Lipschitz trivial, then  $f_t$ is bi-Lipschitz equivalent to $f_{t'}$ for any $t\neq t' \in U$.
\end{theorem}

The above theorem is known as a result of Thom-Levine type and its proof is immediate, since the flow of a Lipschitz vector fields defines a family of bi-Lipschitz homeomorphisms.

\subsection{Weighted homogeneous functions}

Let $w=(w_1,\dots,w_n)$ be an $n$-uple of positive integer numbers. We say that a polynomial function $f(x_1,\dots,x_n)$ is \emph{w-homogeneous} of degree $d$ if $f(sx_1,\dots,sx_n)=s^df(x_1,\dots,x_n)$ for all $s\in\C^*$. Let $H_w^d(n,1)$ denote the space of $w$--homogeneous $n$ variables polynomials of degree $d$.
Let $\mathcal{O}_n$ be the ring of analytic function-germs at the origin $0\in\C^n$ and let $\mathcal{M}_n$ be the maximal ideal of $\mathcal{O}_n$.

\begin{proposition}\label{analytictriviality} Let $F(x_1,\dots,x_n,t)$ be a polynomial function such that: for each $t\in U$, the function $f_t(x_1,\dots,x_n)=F(x_1,\dots,x_n,t)$ is $w$--homogeneous with an isolated singularity at $0\in\C^n$, where $U$ is a domain of $\C$. If, for any $t\in U$, $\Ft$ belongs to the ideal of $\mathcal{O}_n$ generated by $\displaystyle\{x_i\frac{\partial F}{\partial x_j} \ : \ i,j=1,\dots,n\}$ , then $f_{t_1}$ is analytically equivalent to $f_{t_2}$ for any $t_1,t_2\in U$.
\end{proposition}

\begin{proof} Let us denote by $TF$ the ideal of $\mathcal{O}_n$ generated by $\displaystyle\{x_i\frac{\partial F}{\partial x_j} \ : \ i,j=1,\dots,n\}$. It is clear that $H^d_w(n,1)$ can be considered a subset of  the space of $m$--jets $J^m(n,1)$, for $m$ large enough. The set $$A^d_w(n,1)=\{f\in H^d_w(n,1) \ : \ f \ \mbox{has an isolated singularity at origin} \ \}$$ is a Zariski open subset of $H^d_w(n,1)$. In particular, $A^d_w(n,1)$ can be seen as a connected submanifold of the $m$-jets space $J^m(n,1)$. Let $\mathcal{R}(n,n)$ be the group of analytic diffeomorphism-germs $(\C^n,0)\rightarrow(\C^n,0)$. We consider the natural action of $G=j^m(\mathcal{R}(n,n))$ on the manifold $M=J^m(n,1)$ given by $$(j^m(\phi),j^m(f))\mapsto j^m(f\circ\phi).$$ So, given $f\in A$, we get:
\begin{eqnarray*}
\mathrm{codim}(G\cdot f) &=& \mathrm{dim}(M)-\mathrm{dim}(G\cdot f) \\
&=&\mathrm{dim}(M) - \mathrm{dim}(T_f(G\cdot f)) \\
&=& \mathrm{dim}(\frac{\mathcal{M}_n}{\mathcal{M}_n^{m+1}}) - \mathrm{dim}(\frac{TF}{\mathcal{M}_n^{m+1}}) \\
&=& \mathrm{dim}(\frac{\mathcal{M}_n}{TF}) \\
&=& n-1 + \mu(f)
\end{eqnarray*}
where $\mu(f)$ is the Milnor number of $f$. Since the Milnor number of $f\in A$ does not depend on $f$, we get that dimension of $T_f(G\cdot f)$ does not depend on $f\in A$. Let $$P=\{f_t \ : \ t\in U\}.$$ Then,
\begin{enumerate}
\item by assumptions, $T_{f}P\subset T_{f}(G\cdot f)$, for any $f\in P$ and,
\item $\mathrm{dim}(T_f(G\cdot f))$ is constant for $f\in P$, because $P\subset A$.
\end{enumerate}

Now we are ready to apply Mather's Lemma in order to know that $P$ is included in a single $G$--orbit.
\end{proof}

\section{Results}

\begin{proposition} Let $F(x,y,t)$ be a polynomial function such that: for each $t\in U$, the function $f_t(x,y)=F(x,y,t)$ is $w$--homogeneous ($w_2>w_1$) with an isolated singularity at $(0,0)\in\C^2$, where $U$ is an open subset of $\C$. If $f_t$ defines a strongly bi-Lipschitz trivial family of function-germs at origin $(0,0)\in\C^2$, then there exists an algebraic continuous function $k\colon U\rightarrow\C$ such that $$\Ft(x,y,t) - k(t)y\Fy(x,y,t)$$ is identically null on the polar set $\{(x,y,t) \ : \ \Fx(x,y,t)=0\}.$
\end{proposition}

\begin{proof} By assumptions, there exists a Lipschitz vector fields $$v(x,y,t)=(v_1(x,y,t),v_2(x,y,t),1)$$ such that $\displaystyle\frac{\partial F}{\partial v}=0$.
 Let $a_1(t),...,a_r(t)$ be continuous functions such that
 $$\gamma_i(s)=(a_i(t)s^{w_1},s^{w_2},t) \ i=1,\dots,r$$ parameterize the branches of
 $$\Gamma_t=\{(x,y,t) \ : \ \Fx(x,y,t)=0\}.$$ Let
 $$k_i(t)=\frac{\Ft(a_i(t),1,t)}{\Fy(a_i(t),1,t)}.$$ We claim that $k_i(t)=k_j(t)$ for any $i\neq j$. In fact, since $v_2(x,y,t)$ is a Lipschitz function, we have that
 $$|v_2(a_i(t)s^{w_1},s^{w_2})-v_2(a_j(t)s^{w_1},s^{w_2})|\lesssim |s|^{w_1} \ \mbox{as} \ s\to 0.$$ On the other hand,

\begin{eqnarray*}
|v_2(a_i(t)s^{w_1},s^{w_2},t)-v_2(a_j(t)s^{w_1},s^{w_2},t)|&=&|\frac{\Ft(a_i(t)s^{w_1},s^{w_2},t)}
{\Fy(a_i(t)s^{w_1},s^{w_2},t)}- \frac{\Ft(a_j(t)s^{w_1},s^{w_2})}{\Fy(a_j(t)s^{w_1},s^{w_2})}| \\
&=& |k_i(t)-k_j(t)||s|^{w_2}.
\end{eqnarray*}
Since $w_1>w_2$, we have $k_i(t)=k_j(t).$

Let us denote $k(t)=k_1(t)=\cdots=k_r(t)$. Now, let us fix $t$. We see that the function
 $$\Ft(x,y,t) - k(t)y\Fy(x,y,t)$$ is identically null on each branch of $\Gamma_t$ at $(0,0,t)$.
\end{proof}

 \begin{theorem}\label{main} Let $F(x,y,t)$ be a polynomial function such that: for each $t\in U$, the function $f_t(x,y)=F(x,y,t)$ is $w$--homogeneous ($w_2>w_1$) with an isolated singularity at $(0,0)\in\C^2$, where $U\subset\C$ is a domain. If $\{f_t \ : \ t\in U \}$, as a family of function-germs at $(0,0)\in\C^2$, is strongly bi-Lipschitz trivial , then $f_{t_1}$ is analytically equivalent to $f_{t_2}$ for any $t_1,t_2\in U$.
 \end{theorem}

 \begin{proof}  Let $k(t)$ be given by the above proposition, hence we have an algebraic continuous function $k\colon U\rightarrow\C$ such that $$\Ft(x,y,t) - k(t)y\Fy(x,y,t)$$ is identically null on the polar set $\{(x,y,t) \ : \ \Fx(x,y,t)=0\}.$*

 Let $t_0\in U$. Using the Newton-Puiseux Parametrization Theorem, there exist two open $0\in V\subset\C$ and $t_0\in U'\subset U$ such that
 $$s\mapsto s^N+t_0$$ maps $V$ onto $U'$ and the function $$s\mapsto k(s^N+t_0)$$ is analytic for some $N\in\N$. Let us denote $G(x,y,s)=F(x,y,s^N+t_0)$, $g_s(x,y)=G(x,y,s)$ and $\tilde{k}(s)=k(s^n+t_0)$. Thus, we have an analytic funtion $\tilde{k}\colon V\rightarrow\C$ such that $$\Gs(x,y,s) - Ns^{n-1}\tilde{k}(s)y\Gy(x,y,s)$$ is identically null on the polar set $\{(x,y,s) \ : \ \Gx(x,y,s)=0\}.$ Let $$P_1(x,y,s),\dots,P_r(x,y,s)$$ be such that
 they define the analytic irreducible factors of $\Gx(x,y,s)$ in $\mathcal{O}_3$. Let $j(s)=Ns^{n-1}\tilde{k}(s)$. Let
 $$\alpha_i=\max\{\alpha\in\N \ : P_i^{\alpha_i} \ \mbox{divides} \ \Gx \ \mbox{in} \ \mathcal{O}_3\}.$$
 By hypothesis, we have an integer number $\beta_i\geq 1$ such that
 $$\Gs(x,y,s)=u(x,y,s)P_i^{\beta_i}(x,y,s)+j(s)y\Fy(x,y,t)$$ with $u=u(x,y,s)\in\mathcal{O}_3$. Moreover, we can suppose that $P_i^{\beta_i}$ does not divide $u$ in $\mathcal{O}_3$. We should show that $\beta_i\geq\alpha_i$. If $\alpha_i= 1$, we have nothing to do. Thus, let us consider $\alpha_i>1$. It follows from
 $$\Gs(x,y,s)=u(x,y,s)P_i^{\beta_i}(x,y,s)+j(s)y\Gy(x,y,s)$$ that
 \begin{equation}\label{eq1}
 \frac{\partial^2G}{\partial x \partial s}=\frac{\partial u}{\partial x}P_i^{\beta_i}+\beta_iuP_i^{\beta_i-1}\frac{\partial P_i}{\partial x}+ jy\frac{\partial^2G}{\partial x\partial y}.
 \end{equation}
 Now, since $P_i^{\alpha_i}$ divides $\displaystyle\Gx$, we have that $P_i^{\alpha_i-1}$ divides $\displaystyle\frac{\partial^2G}{\partial s \partial x}$ and $\displaystyle\frac{\partial^2G}{\partial y \partial x}$, hence, by using eq. \ref{eq1}, we get that
 \begin{equation}\label{eq2}
 P_i^{\alpha_i-1} \ \mbox{divides} \ P^{\beta_i-1}(\frac{\partial u}{\partial x}P_i+\beta_iu\frac{\partial P_i}{\partial x}).
 \end{equation}
 Since $P_i$ does not divide neither $u$ and $\displaystyle\frac{\partial P_i}{\partial x}$, it follows from eq. \ref{eq2} that $\beta_i\geq\alpha_i$.

 Once $\beta_i\geq \alpha_i$, we have that $\Gs$ belongs the ideal of $\mathcal{O}_3$ generated by $$\{x\Gx,x\Gy,y\Gx,y\Gy\}.$$ Thus, it comes from Mather's Lemma or the analytic version of Thom-Levine result that $g_{s_1}$ is analytically equivalent to ${g_{s_2}}$ for any $s_1,s_2\in V$. It means that, $f_{t_1}$ is analytically equivalent to $f_{t_2}$ for any $t_1,t_2\in U'$. Finally, it follows from the connectivity of $U$  that $f_{t_1}$ is analytically equivalent to $f_{t_2}$ for any $t_1,t_2\in U$.

 \end{proof}

 In the above proof, from where is marked *, if $\Fx$ were analytically reduced, this theorem could be proved in the following way. We fix $t\in U$ and, since $\displaystyle\Fx(x,y,t)$ is analytically reduced, there exists an analytic function $u(x,y)$ such that
 $$\Ft(x,y,t)-k(t)y\Fy(x,y,t)=u(x,y)\Fx(x,y,t).$$ It comes from this equation that $u(x,y)$ is $w$--homogeneous of degree $w_1$, in particular $u(0,0)=0$. Thus, $\Ft$ belongs to the ideal of $\mathcal{O}_2$ generated by
 $$\{x\Fx(x,y,t),x\Fy(x,y,t),y\Fx(x,y,t),y\Fy(x,y,t)\}.$$ Finally, since $f_t(x,y)=F(x,y,t)$ is $w$--homogeneous for all $t\in U$, it follows from Proposition \ref{analytictriviality} that $\{f_t \ : \ t\in U\}$ defines a family of function-germs at origin $(0,0)\in\C^2$ such that $f_{t_1}$ is analytically equivalent to $f_{t_2}$ for any $t_1,t_2\in U$.

The above argument allows us to extend the Theorem \ref{main} to  $w$-homogeneous polynomials in $n$ variables ($w_1>\dots>w_n$) with the additional hypothesis that the ideal $$\displaystyle\{\frac{\partial F}{\partial x_1},\dots,\frac{\partial F}{\partial x_{n-1}}\}$$ is radical. However, the following example shows that, for $n\geq 3$ variables, if we remove some hypothesis above this result is not true at all.

\begin{example} $f_t(x,y,z)=x^4+y^4+z^k+tx^2y^2$ is strongly bi-Lipschitz trivial (for $t$ close to $0$), but $f_{t_1}$ is not analytically equivalent to $f_{t_2}$ when  $t_1\neq t_2$.
\end{example}

\end{document}